\documentclass[article,12pt]{article}

\usepackage[english]{babel}

\usepackage[a4paper, total={6in, 9in}]{geometry}

\usepackage[figurename=Fig.]{caption}

\usepackage{amssymb,amsfonts,amsmath,amsthm}
\usepackage{graphicx}
\usepackage[colorlinks=true, allcolors=black]{hyperref}

\usepackage{graphicx}
\usepackage{amssymb,amsmath}
\usepackage{amsfonts}
\usepackage{latexsym}
\usepackage{url}
\usepackage{color}
\usepackage{multicol}
\usepackage{booktabs}

\def\R{\mathbb{R}}
\def\Z{\mathbb{Z}}

\def\Mcp{\mathcal{P}}

\newtheorem{theorem}{Theorem}
\newtheorem{corollary}[theorem]{Corollary}
\newtheorem{proposition}[theorem]{Proposition}

\newcounter{example}[section]
\newenvironment{example}[1][]{\refstepcounter{example}\par\medskip
   \noindent {\it Example~\theexample. #1} \rmfamily}{\medskip}

\newenvironment{keywords}{\begin{quote}\small \em 
{\bf Keywords\/}:}{\end{quote}}

\usepackage{authblk}
\title{Comparing perspective reformulations for piecewise-convex optimization}
\author[1]{Renan Spencer Trindade\thanks{Corresponding author

Email addresses:
\href{mailto:rst@lix.polytechnique.fr}{rst@lix.polytechnique.fr} (Renan Spencer Trindade),

\noindent
\href{mailto:dambrosio@lix.polytechnique.fr}{dambrosio@lix.polytechnique.fr} (Claudia D'Ambrosio),
\href{mailto:frangio@di.unipi.it}{frangio@di.unipi.it} (Antonio Frangioni),

\noindent
\href{mailto:gentile@iasi.cnr.it}{gentile@iasi.cnr.it} (Claudio Gentile).
}}
\author[1]{Claudia D'Ambrosio}
\author[2]{Antonio Frangioni}
\author[3]{Claudio Gentile}

\affil[1]{LIX, CNRS, École Polytechnique, Institut Polytechnique de Paris

Route de Saclay, Palaiseau, 91128, France
}

\affil[2]{Dipartimento di Informatica, Università di Pisa

Largo B. Pontecorvo 3, Pisa, 56127, Italy
}

\affil[3]{Istituto di Analisi dei Sistemi ed Informatica ``Antonio Ruberti'', C.N.R.

Via dei Taurini 19, Rome, 00185, Italy
}

\date{}

\begin{document}
\maketitle

\begin{abstract}
Our study is motivated by the solution of Mixed-Integer Non-Linear Programming (MINLP) problems with separable non-convex functions via the Sequential Convex MINLP technique, an iterative method whose main characteristic is that of solving, for bounding purposes, piecewise-convex MINLP relaxations obtained by identifying the intervals in which each univariate function is convex or concave and then relaxing the concave parts with piecewise-linear relaxations of increasing precision. This process requires the introduction of new binary variables for the activation of the intervals where the functions are defined. In this paper we compare the three different standard formulations for the lower bounding subproblems and we show, both theoretically and computationally, that---unlike in the piecewise-linear case---they are not equivalent when the perspective reformulation is applied to reinforce the formulation in the segments where the original functions are convex.
\end{abstract}

\section{Introduction}

The Sequential Convex Mixed-Integer Non-Linear Programming  (SC-MINLP) technique, introduced in \cite{DLW09,DAmbrosio2008}, aims at solving non-convex Mixed-Integer Non-Linear Programs (MINLPs) in which non-convexity is restricted to sums of univariate functions, i.e.,
%
\begin{align*}
 \min \; &\textstyle \sum_{j \in N} c_jx_j \\
 & \textstyle
   f_i(x) + {\sum_{j \in H(i)} g_{ij}(x_j) \leq 0 } & i \in M \\
 & l_j \leq x_j \leq u_j &  j \in N \\
 & x_j \in \Z & j \in I
\end{align*}
where the multivariate functions $f_i : \R^n \to \R$ are convex while the univariate functions $g_{ij} : \R \to \R$ are non-convex. All sets $M$, $N$, $I \subseteq N$, and $H(i) \subseteq N$ are finite, and $l_j$ and $u_j$ are finite bounds for the $x_j$ with $j \in H(i)$, i.e., those appearing in some $g_{ij}$ function for some $i \in M$. SC-MINLP is less generic than solvers like Baron \cite{ts:05,sahinidis:baron:21.1.13}, Couenne \cite{Belotti2009}, SCIP \cite{BestuzhevaEtal2021ZR}, but, as it exploits some mathematical properties of the class of MINLPs at hand, it can be more efficient on the problems it can handle.

In the SC-MINLP technique, a lower bound is obtained by solving a convex MINLP obtained as follows.  Let us consider a pair $i,j$ where $i \in M, j \in H(i)$. We first compute the $s(ij)+1$   breakpoints $l_ j = l^1_{ij} < l^2_{ij} < ... < l^{s(ij)}_{ij} < l^{s(ij)+1}_{ij}=u_j$ that separate subintervals of the domain of $x_{j}$ where $g_{ij}$ is concave or convex, using the second derivative of $g_{ij} (x_j)$. Then, we define $S(ij)=\{1,\ldots,s(ij)\}$,  $\check{S}(ij)=\{s \in S(ij) ~|~ g_{ij} \mbox{ is convex in } [l^s_{ij}, l^{s+1}_{ij}] \}$, and $\hat{S}(ij)=\{s \in S(ij) ~|~ g_{ij} \mbox{ is concave in } [l^s_{ij}, l^{s+1}_{ij}] \}$. In the intervals in $\hat{S}(ij)$, we replace the concave function $g_{ij}$ by its convex envelope (a linear function), thus obtaining a piecewise-convex MINLP whose continuous relaxation is convex. In the original description of SC-MINLP \cite{DLW09,DAmbrosio2008,DAmbrosio2019}, the \emph{incremental model}, described in the next section, was used to formulate the piecewise-convex functions. Motivated by this setting we study and compare different formulations of piecewise-convex MINLP problems. This has been done in the piecewise-\emph{linear} case in \cite{Croxton2003}, showing that the three ``textbook'' formulations, namely the incremental, the convex combination, and the multiple choice formulation are equivalent in terms of continuous relaxations bounds. In \cite{VAN2010}, new so-called logarithm formulations are introduced, one for each textbook formulation, that use  less binary variables so as to potentially speed-up the Mixed Integer Linear Programming (MILP) solvers. In this paper we prove that, somewhat surprisingly, the results for the linear case do not carry over to the nonlinear one, in that one formulation---incidentally, the one employed in \cite{DLW09,DAmbrosio2008,DAmbrosio2019}---is \emph{weaker} than the other two even when all are strengthened by the use of the \emph{Perspective Reformulation} (PR) technique \cite{FrGe06a,FrFG16}. Our computational experience proves that the stronger formulation of the subproblem directly translates in better performances of the SC-MINLP technique.

\section{MIP Models for SC-MINLP}
We now discuss the generalization of the textbook formulations for piecewise-linear functions to piecewise-convex ones. 

\subsection{Incremental model}
In the \emph{Incremental Model} (IM), which was used in the original version of SC-MINLP \cite{DLW09,DAmbrosio2008,DAmbrosio2019}), each sub-interval $[l_{ij}^s, l_{ij}^{s+1}]$ has a segment load variable $x^s_{ij}$ which assumes value zero unless $x^{s-1}_{ij}$ reaches its maximum value, that is, $x_{ij}^{s} > 0$ only if $x_{ij}^{s-1} = l_{ij}^{s} - l_{ij}^{s-1}$. The binary variables $y_{ij}^{s}$ are defined by the condition that $y_{ij}^{s}=1$ if $x_{ij}^{s}>0$, and $y_{ij}^{s}=0$ otherwise. The ``plain'' version of the formulation would then have terms $g_{ij}(l^s_{ij} + x^s_{ij}) - y^s_{ij}  g_{ij}(l^s_{ij})$ in each constraint $i \in M$ for $j \in H(i)$ and $s \in \check{S}(ij)$. However, it is well-known that the \emph{semi-continuous} variable $x_{ij}^{s}$, governed by the binary variables $y_{ij}^{s}$, lends itself to the PR technique \cite{DAmbrosio2019} that can considerably strengthen the continuous relaxation of the problem. This comes at the cost of ``more nonlinear'' functions and possibly the introduction of further auxiliary variables, resulting, e.g., in
\begin{align} 
 \min \; & \textstyle \sum_{j \in N} c_jx_j \label{eq:Inc-obj} \\
 & \textstyle
   \bar{f}_i(x) + {\sum_{j \in H(i)} \sum_{s \in \check{S}(ij)} z^s_{ij} \leq 0 }
 & i \in M \label{eq:Inc-const} \\
 %
 %
 & z^s_{ij} \geq [g_{ij} (l^s_{ij} + x^s_{ij}/y^s_{ij}) - g_{ij}(l^s_{ij})] y^s_{ij} 
 & s \in \check{S}(ij) \;,\; j \in H(i) \;,\; i \in M
   \label{eq:Inc-const2} \\
 & y^1_{ij} = 1 \;\;,\;\; y^{s(ij)+1}_{ij} = 0
 & j \in H(i) \;,\; i \in M
   \label{eq:Inc-firstlast} \\
 & \textstyle
   x_j = l_j + \sum_{s \in S(ij)} x^s_{ij} & j \in H(i) \;,\; i \in M \\
 & (l^{s+1}_{ij} \!-\! l^s_{ij} ) y_{ij}^{s+1} \leq x_{ij}^s \leq (l^{s+1}_{ij} \!-\! l^s_{ij}) y_{ij}^s
 & s \in S(ij) \;,\; j \in H(i) \;,\; i \in M
   \label{eq:boundsxijs2} \\
 & y_{ij}^s \in \{0, 1\} & s \in S(ij) ,\; j \in H(i) ,\; i \in M \label{eq:Inc-varY} \\
 & x_j \in \Z & j \in I
 \label{eq:Inc-int}
\end{align}
where $\bar{f}_i = f_i(x) + \sum_{j \in H(i)} g_{ij}(l^1_{ij}) + \sum_{s \in \hat{S}(ij)} \alpha^s_{ij} x^s_{ij}$ and $\alpha^s_{ij} = (g_{ij}(l_{ij}^{s+1})-g_{ij}(l_{ij}^s))/(l_{ij}^{s+1} - l_{ij}^s)$ is the slope of the linear function for ``concave'' intervals $s \in \hat{S}(ij)$. In particular, \eqref{eq:Inc-const2} implements the PR of the ``plain'' constraint. In \eqref{eq:Inc-firstlast}, $y^{s(ij)+1}_{ij} = 0$ is just a syntactic trick to write \eqref{eq:boundsxijs2} in an uniform way avoiding the ``border effect'' of the last segment. Instead, $y^1_{ij} = 1$ is significant: the constraint may be avoided, allowing $y^1_{ij} = 0$ to happen when $x_j = l_j$, but, by dint of having less solutions, the continuous relaxation is clearly stronger if we fix the variable.

\subsection{Multiple-choice model}

An alternative formulation is the \emph{Multiple-Choice Model} (MCM) where, for each $s \in S(ij)$, the load variable $x^s_{ij}$ defines the total load $x^s_{ij}=x_j$ and $y^{s}_{ij}=1$ if $x_j$ lies on the sub-interval $[l_{ij}^s, l_{ij}^{s+1}]$, and $x^s_{ij} = y^{s}_{ij} = 0$ otherwise. In this formulation, exactly one $y^{s}_{ij}$ will equal one. In this case, the ``plain'' terms in the constraints would have the form $g_{ij} (x^s_{ij}) - y^{s}_{ij} g_{ij}(0)$: strengthened with the PR technique yields
\begin{align} \hspace{0.5cm}
 & z^s_{ij} \geq [g_{ij} (x^s_{ij}/y^s_{ij}) - g_{ij}(0)]y^s_{ij} 
 & s \in \check{S}(ij) \;,\; j \in H(i) ;,\; i \in M
   \label{eq:MC-const2} \\
 & \textstyle
   x_j = \sum_{s \in S(ij)} x^s_{ij} & j  \in H(i) \;,\; i \in M
   \label{eq:redefxj2} \\
 & l^{s}_{ij}\! y_{ij}^{s} \leq x_{ij}^s \leq l^{s+1}_{ij}\! y_{ij}^s
 & s \in S(ij) \;,\; j \in H(i) \;,\; i \in M \\
 & \textstyle
   \sum_{s \in S(ij)} y_{ij}^{s} = 1
 & i \in M \;,\; j \in H(i) \label{eq:Inc-SumY}
\end{align}
together with \eqref{eq:Inc-obj}, \eqref{eq:Inc-const}, \eqref{eq:Inc-varY} and \eqref{eq:Inc-int}, except that now $\bar{f}_i$ in \eqref{eq:Inc-const} is rather $f_i(x) + \sum_{j \in H(i)} g_{ij}(0) \sum_{s \in \check{S}(ij)} y_{ij}^s + \sum_{s \in \hat{S}(ij)} (\alpha^s_{ij} x^s_{ij} + (g_{ij}(l_{ij}^s)-\alpha^s_{ij}l^{s}_{ij} ) y_{ij}^s)$.

\smallskip
\noindent
The third ``textbook'' formulation is the \emph{Convex Combination Model} (CCM) \cite{Croxton2003}, which can be succinctly described as follows. On each $s \in S(ij)$, $x^s_{ij}$ from MCM is replaced by the pair of (convex combination) variables $\lambda^s_{ij}$ and $\mu^s_{ij}$, i.e., $x^s_{ij} = l^s \mu^s_{ij} + l^{s+1} \lambda^{s}_{ij}$, under the constraint $y^s_{ij} = \mu^s_{ij} + \lambda^s_{ij}$. For each continuous solution $(x,y)$ of the MCM, it is easy to obtain a solution $(\mu,\lambda,y)$ of the CCM with the same value of the constraint \eqref{eq:MC-const2} by just taking (dropping the $ij$ index for simplicity) $\mu^s = \frac{l^{s+1}-x^s/y^s}{l^{s+1}-l^s} y^s=\frac{l^{s+1}y^s-x^s}{l^{s+1}-l^s}$ and $\lambda^s = \frac{x^s/y^s - l^s}{l^{s+1} - l^s} y^s=\frac{x^s - l^sy^s}{l^{s+1} - l^s}$. Since the opposite result is trivial, CCM and MCM have equivalent continuous relaxation bound and we can restrict ourselves to comparing IM with MCM.

\section{Incremental Model vs Multiple Choice Model}

While in the piecewise-linear case IM and CCM have equivalent continuous relaxations~\cite{Croxton2003}, this is no longer true in the general piecewise-convex case. In fact, while MCM describes the convex envelope of the piecewise function, IM does not. The first point comes from \cite[\S 3.1]{FGH20}, which considers a $x \in \R^n$ partitioned into $x = [x^k]_{k \in K}$, where each $x^k$ is either $0$ or belongs to a compact convex set (for our purposes, a polytope $\Mcp^k = \{ \; x^k \;:\; A^k x^k \leq b^k \; \}$, although this is not necessary in general). Let $g^k$ be a convex objective function defined on $\Mcp^k$ with $g^k(0)=0$ and $c^k$ be a fixed activation cost. For the \emph{alternatives function}
\begin{equation*}
 g(x) =
 \left\{\begin{array}{ll}
         g^k(x^k) + c^k & \mbox{if } x^k \in \Mcp^k \mbox{ and }
                          x^h = 0 \; \forall\; h \in K \setminus \{ k \}
                          \\
         0              & \mbox{if } x = 0 \\
         + \infty       & \mbox{otherwise}
        \end{array}\right.
 \;\; .
 \label{eq:fgen}
\end{equation*}
the convex envelope of $g$ can be described as
\begin{align*}
 \overline{co} \, g(x) \,=\,
 \min \; & \textstyle
           \sum_{k \in K} \theta^k g^k(x^k / \theta^k) \\
 & \textstyle \sum_{k \in K} \theta^k \leq 1 \\
 & A^k x^k \leq b^k \theta^k \;,\; \theta_k \geq 0
 & \quad k \in K
\end{align*}
This construction is exactly the one applied to each function $g_{ij}$ when formulating the Multiple Choice Model; as a consequence

\begin{corollary} \label{cor:IMvsMCM}
In the Multiple Choice Model, constraints \eqref{eq:MC-const2}--\eqref{eq:Inc-SumY} describe the convex envelope of each function $g_{ij}$.
\end{corollary}

Thus, the lower bound provided by MCM is not lower than the lower bound produced by IM. The following example shows that the MCM bound can be strictly better than the IM bound.

\begin{example}
\label{ex:convex-concave}
Consider the problem $\min\{ \, y \,:\, y \geq g(x) \;,\; 0 \leq x \leq 2\pi \, \}$ for $g(x) = -\sin(x)$: the feasible region $[0,2\pi]$ is divided in the two intervals [0, $\pi$], where $g$ is convex, and [$\pi$, $2\pi$], where $g$ is concave, as shown in Figure \ref{fig:ex2_1}. The original function $g(x)$ is represented by the black dotted line, while the function after the piecewise-convex reformulations is shown in a solid line, as well as its feasible region in gray.

\begin{figure}[ht]
\centering
\begin{minipage}{.49\textwidth}
\includegraphics[width=\textwidth]{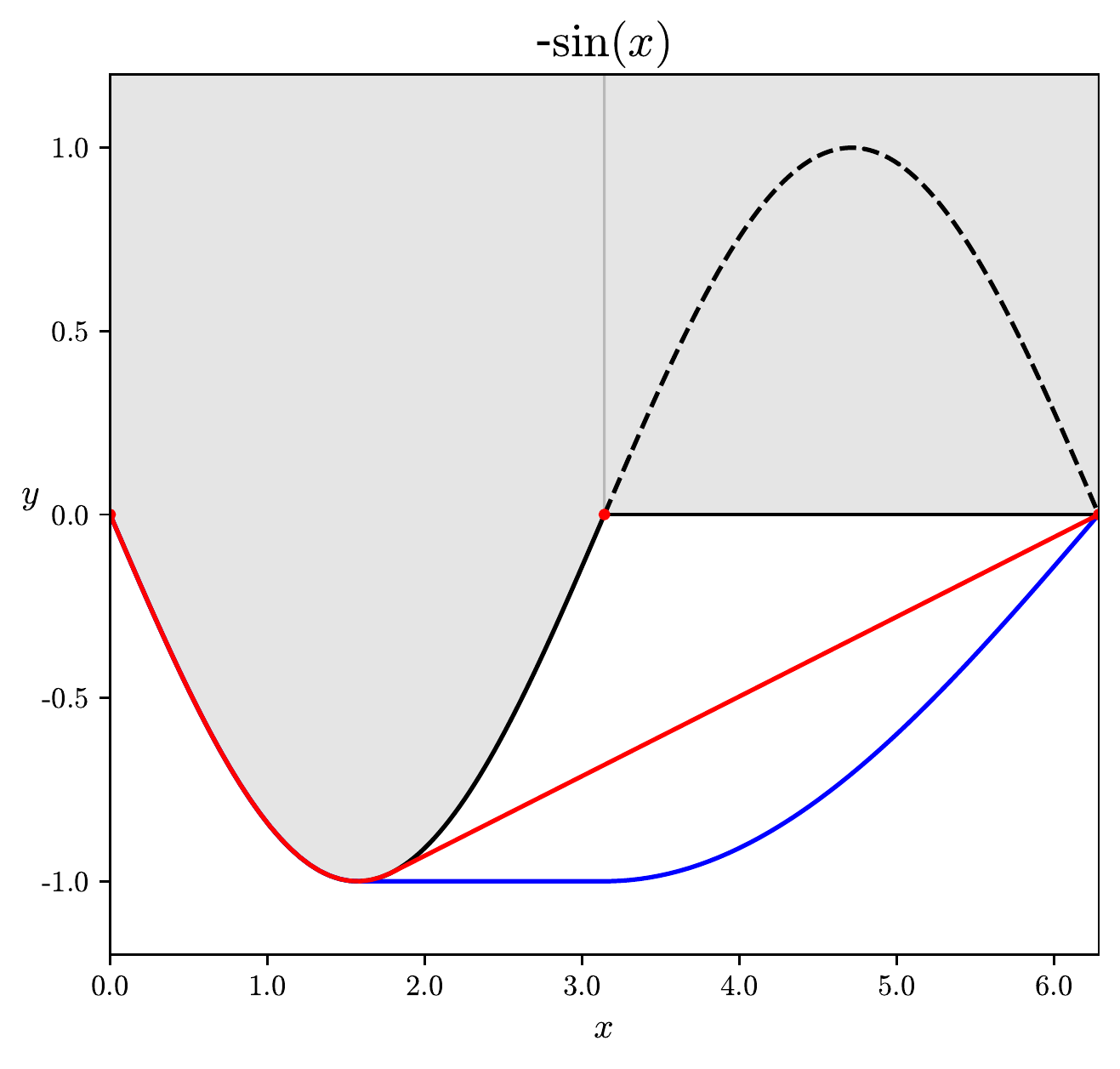}
\caption{Example where the continuous relaxations of MCM and IC are different.}
\label{fig:ex2_1}
\end{minipage}%
\hspace{.01\textwidth}
\begin{minipage}{.49\textwidth}
\includegraphics[width=\textwidth]{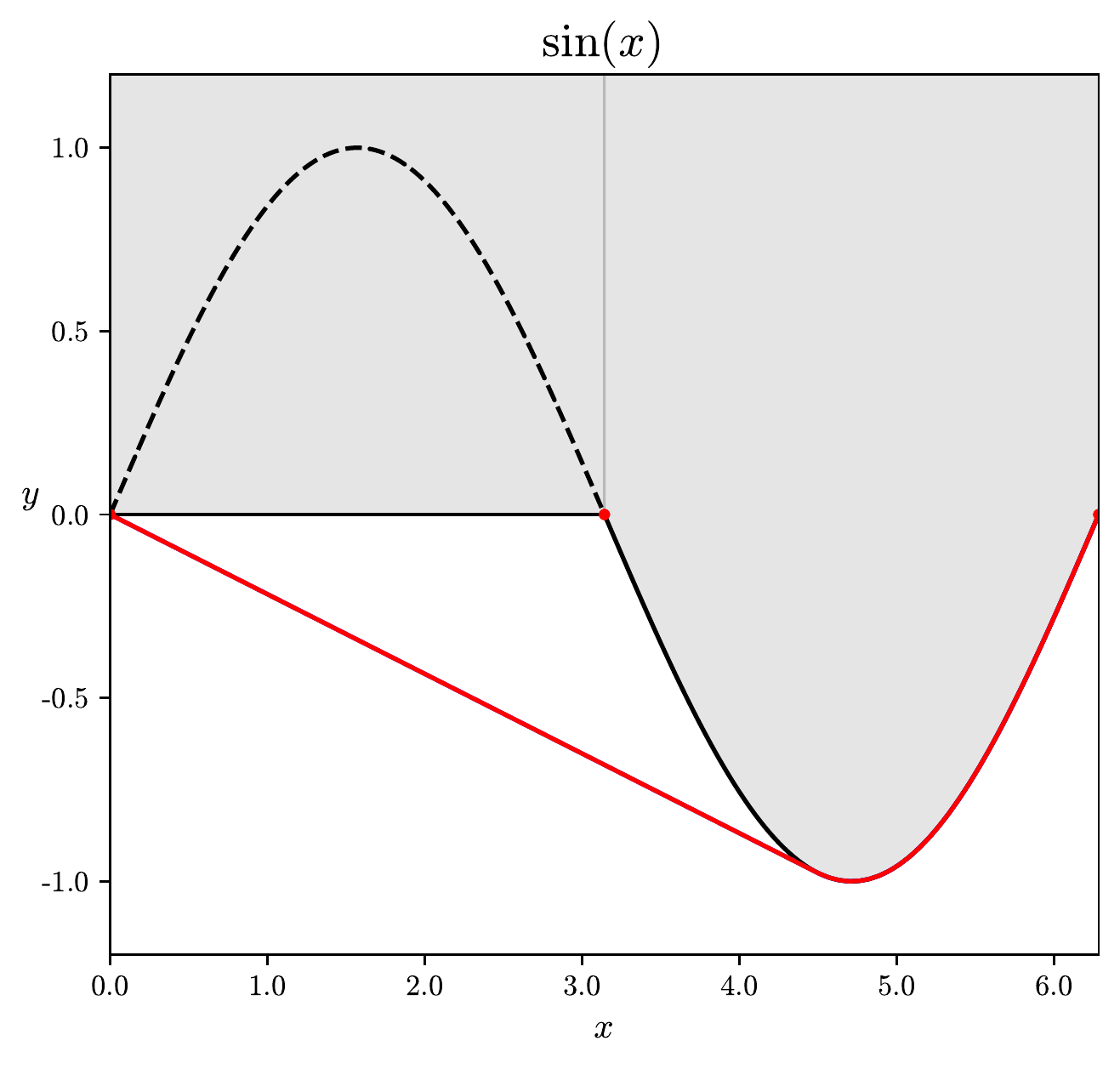}
\caption{Example where the continuous relaxations of MCM and IC are equivalent.}
\label{fig:ex2_2}
\end{minipage}
\end{figure}

The feasible region to this problem, when integrality constraints are present, is the same for both MCM and IC: however, this is no longer true for the continuous relaxation. To illustrate this, we fix the value of the variable $x$ for each point in the interval and plot the optimal value of the continuous relaxation of IM and MCM, respectively, in blue and red. Thus, the example proves that the continuous relaxation of IM can be weaker than that of MCM. However, in a particular case the two formulations can be proven to be equivalent:

\end{example}

\begin{proposition}
\label{prop:two-segments}
If the domain of the $g_{ij}$ can be partitioned into two subsets $[l^1, l^2]$ and $[l^2,l^3]$ so that $g_{ij}$ is concave in the former and convex in the latter, then the continuous relaxations of MCM and IM are equivalent.
\end{proposition}

This is depicted in Figure \ref{fig:ex2_2} for $g(x) = sin(x)$. The proof of this proposition is not difficult but rather long and technical, and therefore we postpone it to the appendix.

\section{Computational Results}
We test our approach on two classes of MINLPs problems with the required structure: the Non-linear Continuous Knapsack (NCK) and the Uncapacitated Facility Location (UFL) problems. To solve the piecewise-convex MINLP we use Perspective Cuts \cite{FrGe06a,DAmbrosio2019}, due to their advantageous property of linearising the ``complex'' terms $[g_{ij} (l^s_{ij} + x^s_{ij}/y^s_{ij}) - g_{ij}(l^s_{ij})] y^s_{ij} $ and $[g_{ij} (x^s_{ij}/y^s_{ij}) - g_{ij}(0)]y^s_{ij}$ in the formulations, as well as their good performances w.r.t.~possible alternatives \cite{FrGe09a}. Since our test problems contain no other nonlinearities save for the $g_{ij}$, this reduces the Convex MINLP in SC-MINLP to a MILP, that we then solve with CPLEX 12.10 with a time limit of 10000 seconds. Our results only consider the MILP of the first iteration of SC-MINLP, in which $\hat{S}(ij)$ are the original ones. Subsequent iterations would split some of the intervals in some $\hat{S}(ij)$, which is unlikely to qualitatively change the results of interest here, which are only related to how the formulation behaves on the $\check{S}(ij)$ intervals (that do not change along iterations). 

\subsection{Non-Linear Continuous Knapsack}

The NCK problem, already considered in \cite{DAmbrosio2019}, is
\begin{equation*}
 \textstyle
 \max\big\{ \, \sum_{j \in N} p_j \,:\, \sum_{j \in N} w_j x_j \leq C \;\;,\;\; p_j \leq g_j( x_j ) \;,
            \; 0 \leq x_j \leq U_j \quad j \in N \, \big\}
\end{equation*}
We report results on two sets of instances. In all the cases we randomly generated 10 instances for each value of $|N| \in \{10, 20, 50, 100, 200, 500, 1000\}$, the weight $w_j$ were uniformly randomly drawn from $[1,100]$, the capacity is $C = 50\sum_j w_j$, and $U_j = 100$ for all $j$. Then, in the first group
\[
 g_j( x_j ) = \frac{c_j}{1+b_j \exp (-a_j (x_j + d_j))}
\]
where, for all $j \in N$, $a_j \in [0.1,0.2]$,  $b_j \in [0,100]$, $c_j \in [0,100]$, and $d_j \in [-100,0]$ were uniformly drawn in the corresponding intervals. This results in at most 2 intervals for each $g_j$. In the second set, instead,
\[
 g_j( x_j ) = 7.5\sin\left(\pi\dfrac{x_j-10}{40}\right) - 15\cos\left(\pi\dfrac{x_j-10}{80} \right) + 19.5
\]
which rather yields 4 intervals.

\smallskip
\noindent
In Table \ref{Table1}, each row shows the average results of the 10 instances corresponding to one value of $|N|$; in particular, for each formulation we report the objective function of the best solution found, the total execution time, and the total number of Perspective Cuts (PC) generated. We then concentrate on the root node relaxation (``Relax'' columns) for which we report again the execution time and the number of PCs, together with the relative gap of the correspondig bound w.r.t.~the optimal integer solution.

\begin{table}[ht!]
\centering
\caption{Computational results for Non-linear Continuous Knapsack problem}
\label{Table1}
\setlength{\tabcolsep}{3.2pt}
\scalebox{1}{
    \begin{tabular*}{1.01\linewidth}{rr|rr|rr|rrr|rrr}
    \hline
    \multicolumn{2}{c}{\textsc{ Inst.}} &
    \multicolumn{2}{c}{\textsc{ IM}} &
    \multicolumn{2}{c}{\textsc{ MCM}} &
    \multicolumn{3}{c}{\textsc{ IM relax.}} &
    \multicolumn{3}{c}{\textsc{ MCM relax.}} \cr
    \cmidrule(l){1-2} 
    \cmidrule(l){3-4} 
    \cmidrule(l){5-6} 
    \cmidrule(l){7-9} 
    \cmidrule(l){10-12} 
    \it Int. & \it Size & \it Time & \it Cuts & \it Time & \it Cuts & \multicolumn{1}{c}{\it Gap} & \it Time & \it Cuts & \multicolumn{1}{c}{\it Gap} & \it Time & \it Cuts \cr
    \hline
    2 & 10 & 0.02 & 114.70 & 0.03 & 105.60 & 0.48 & 0.01 & 50.30 & 0.48 & 0.01 & 50.30 \\
    2 & 20 & 0.03 & 187.40  & 0.03 & 179.80 & 0.18 & 0.01 & 92.20 & 0.18 & 0.02 & 92.20 \\
    2 & 50 & 0.05 & 448.20 & 0.05 & 448.20 & 0.02 & 0.02 & 246.10 & 0.02 & 0.02 & 246.10 \\
    2 & 100 & 0.09 & 759.00 & 0.09 & 759.50 & 0.00 & 0.04 & 499.00 & 0.00 & 0.05 & 499.00 \\
    2 & 200 & 0.21 & 1614.50 & 0.22 & 1635.90 & 0.00 & 0.09 & 989.40 & 0.00 & 0.08 & 989.40 \\
    2 & 500 & 0.45 & 3293.90 & 0.45 & 3202.20 & 0.00 & 0.22 & 2504.50 & 0.00 & 0.22 & 2504.50 \\
    2 & 1000 & 1.12 & 5949.60 & 1.00 & 5896.30 & 0.00 & 0.45 & 5039.40 & 0.00 & 0.43 & 5039.40 \\
    \hline
    4 & 10 & 0.06 & 348.40 & 0.04 & 239.70 & 1.31 & 0.01 & 108.10 & 0.17 & 0.01 & 78.80 \\
    4 & 20 & 0.09 & 533.90 & 0.04 & 325.90 & 1.00 & 0.02 & 225.40 & 0.03 & 0.02 & 155.10 \\
    4 & 50 & 0.41 & 1546.10 & 0.16 & 886.70 & 0.80 & 0.04 & 501.90 & 0.01 & 0.04 & 360.30 \\
    4 & 100 & 0.91 & 2332.30 & 0.26 & 1416.30 & 0.83 & 0.10 & 1058.90 & 0.00 & 0.07 & 733.50 \\
    4 & 200 & 3.10 & 4171.30 & 0.54 & 2369.80 & 0.85 & 0.22 & 2255.20 & 0.00 & 0.15 & 1481.00 \\
    4 & 500 & 20.34 & 8931.70 & 2.40 & 5141.40 & 0.80 & 0.69 & 5611.90 & 0.00 & 0.42 & 3613.30 \\
    4 & 1000 & 174.67 & 18249.10 & 4.51 & 8480.70 & 0.83 & 1.92 & 11029.20 & 0.00 & 1.20 & 7363.50 \\
    \hline
    \end{tabular*}
}
\end{table}

As expected, both formulations produce the same optimal integer solutions. When the number of intervals is two, also the continuous relaxations are equivalent, down to exactly the same number of generated PCs; this is clearly a case in which the hypotheses of Proposition \ref{prop:two-segments} are satisfied. 
However, for a larger number of intervals the MCM bound is significantly better, which directly translates in much lower total running times, especially as $|N|$ increases.

\subsection{Non-linear Uncapacitated Facility Location}

\begin{figure}[ht]
\centering
\begin{minipage}{.3\textwidth}
\includegraphics[width=\textwidth]{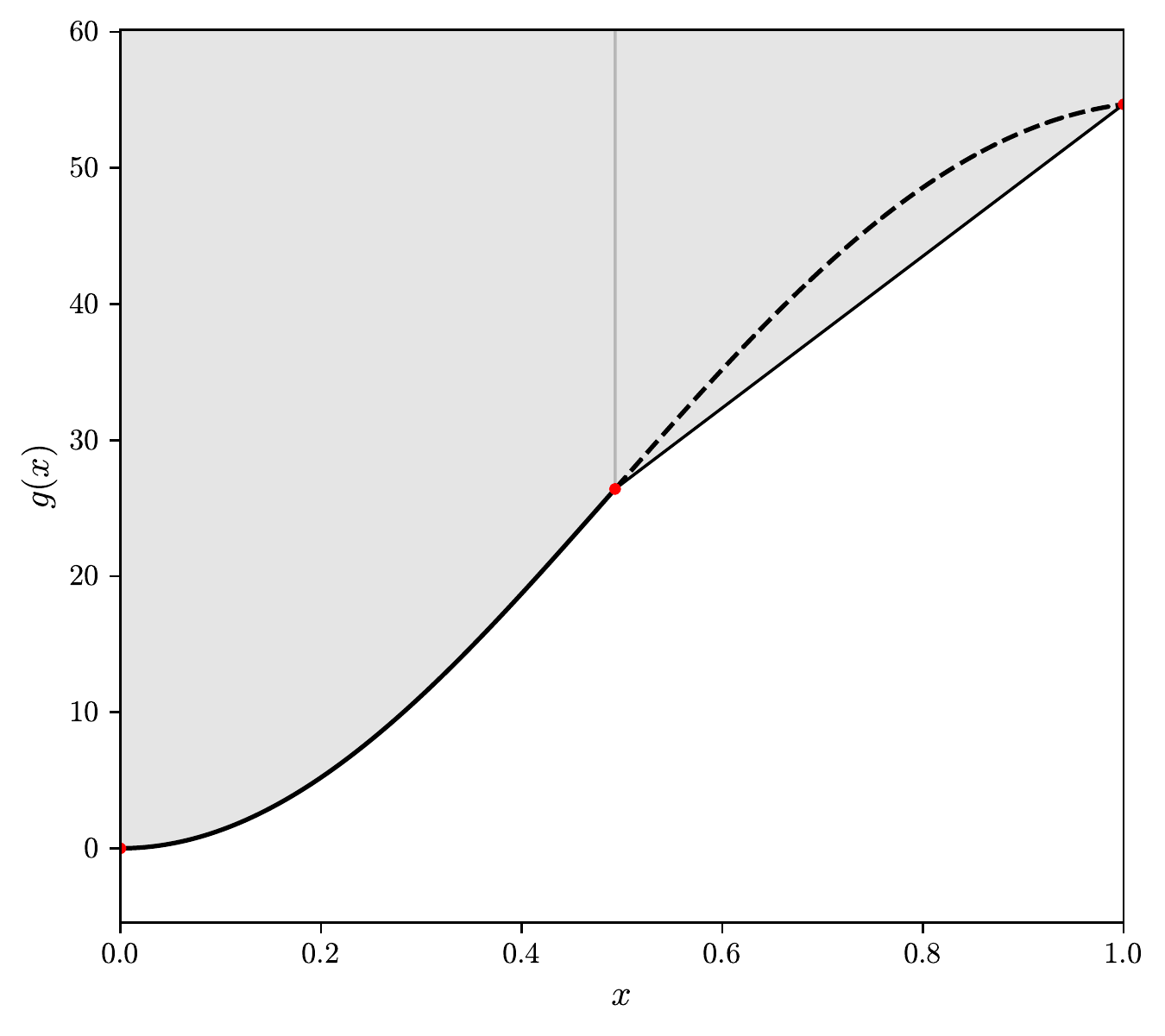}
\caption{Function for UFL instances of type 1 (a=15; b=2; c=1).}
\label{fig:ufl_1}
\end{minipage}%
\hspace{.01\textwidth}
\begin{minipage}{.3\textwidth}
\includegraphics[width=\textwidth]{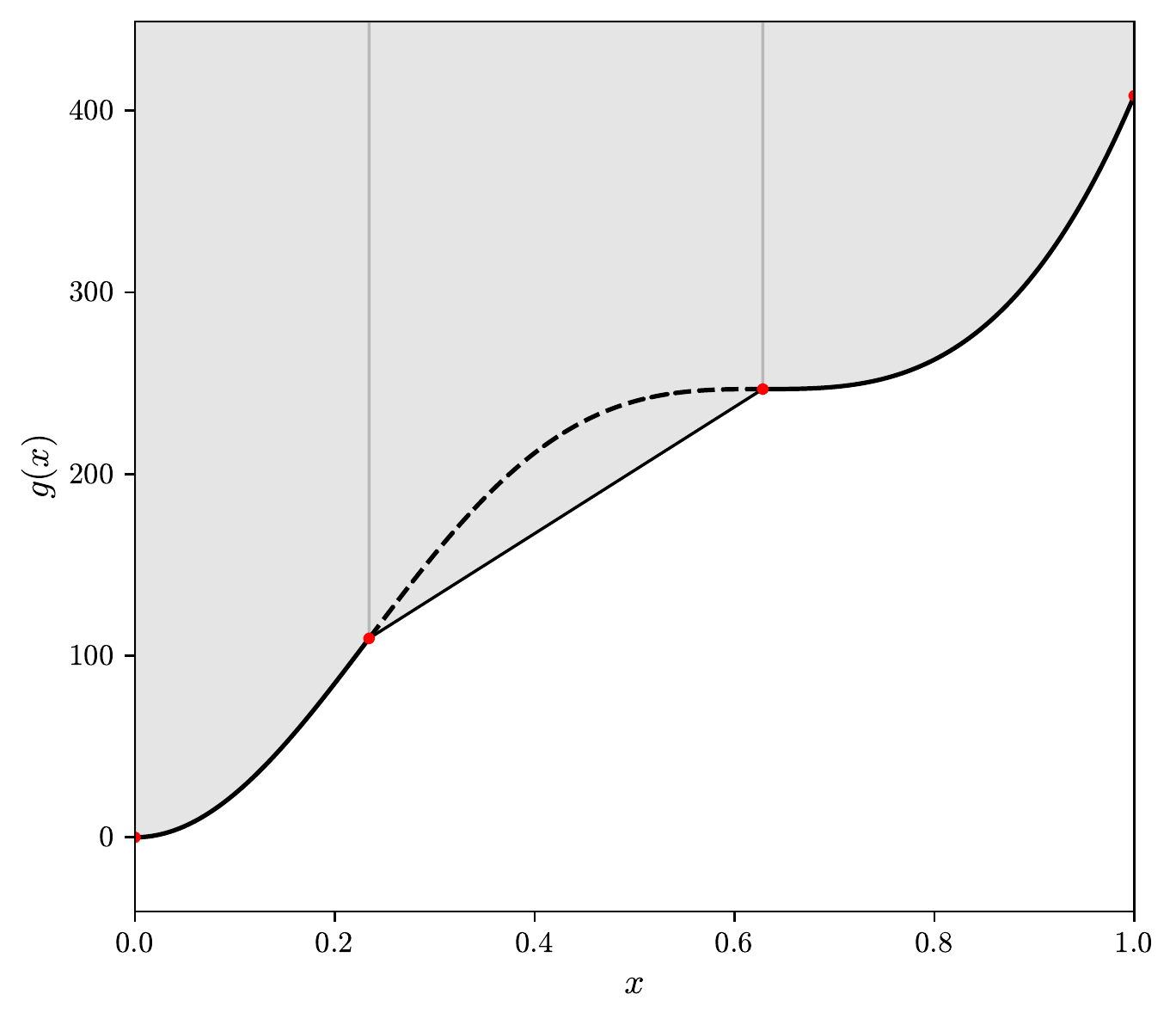}
\caption{Function for UFL instances of type 2 (a=25; b=5; c=5).}
\label{fig:ufl_2}
\end{minipage}
\hspace{.01\textwidth}
\begin{minipage}{.3\textwidth}
\includegraphics[width=\textwidth]{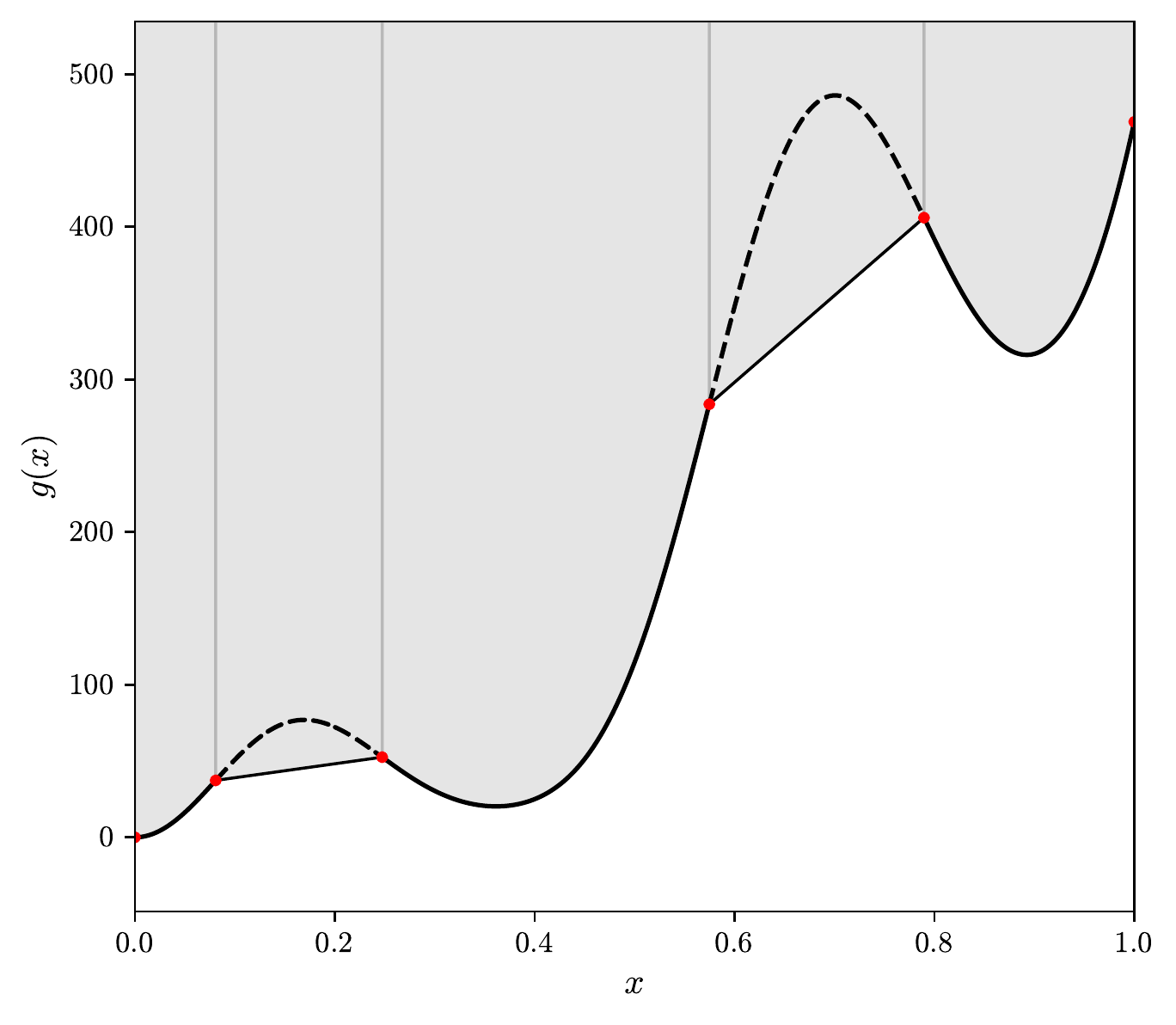}
\caption{Function for UFL instances of type 3 (a=25; b=10; c=5).}
\label{fig:ufl_3}
\end{minipage}
\end{figure}

In the UFL problem, each facility $k \in K$ can satisfy part of the demand required by the consumers in $T$. In this nonlinear version, the shipping cost $s_{kt}$ is defined by a non-convex function, yielding the formulation 
\begin{align} 
 \min & \textstyle \sum_{k \in K} C_k y_k + \sum_{t \in T} \sum_{k \in K} s_{kt} \\
 & s_{kt} \geq a (\sin(b \,\, w_{kt}) + c \,\, w_{kt})^2 & t \in T \,,\, k \in K \label{UFL:cons1} \\
 & \textstyle \sum_{k \in K} w_{kt} = 1 & t \in T \label{UFL:cons2} \\
 & \textstyle 0 \leq w_{kt} \leq y_k & t \in T \,,\, k \in K \label{UFL:cons3} \\
 & \textstyle y_k \in \{0,1\} & k \in K
\end{align}
The variable $w_{kt}$ represents the portion of demand of consumer $t$ satisfied by facility $k$. The binary variable $y_{k}$ activates the facility $k$ at a fixed cost $C_k$ (uniformly generated in [1, 100]). The parameters $a$, $b$ and $c$ are defined for three different types of instances according to Figures \ref{fig:ufl_1}-\ref{fig:ufl_3}, where the functions contain 1, 2, or 3 convex intervals. The number of non-convex constraints is $|K|\cdot|T|$, where we generated three different sizes $(|K|, |T|) \in$ {(6, 12), (12, 24), (24, 48)}. 

Table \ref{Table2} present the computational results for the problem, where each row represents the average of the 10 instances generated for each type. 
In the columns of the integer solutions of the problem, 
the ``\#O'' column reports the number of instances that were solved to optimality of the instances, and ``Gap'' is the value obtained from CPLEX. The other columns in this table have the same meaning as in the previous table.

\begin{table}[!ht]
\caption{Computational results for Non-linear Uncapacitated Facility Location problem}
\label{Table2}
\setlength{\tabcolsep}{3.5pt}
\scalebox{0.78}{
    \begin{tabular*}{1.28\linewidth}{r|rrrr|rrrr|rrr|rrr}
    \hline
     \multicolumn{1}{c}{\textsc{ }} &
    \multicolumn{4}{c}{\textsc{ IM}} &
    \multicolumn{4}{c}{\textsc{ MCM}}  &
    \multicolumn{3}{c}{\textsc{ IM - Relax}} &
    \multicolumn{3}{c}{\textsc{ MCM - relax}}  \\
    \cmidrule(l){2-5} 
    \cmidrule(l){6-9} 
    \cmidrule(l){10-12} 
    \cmidrule(l){13-15}
     \multicolumn{1}{c}{\textsc{\it Inst.}} & 
    \it Time & Gap &\it Cuts & \it \#O & 
    \it Time & Gap &\it Cuts & \it \#O &
    \it Time & \multicolumn{1}{c}{\it Gap } &  \it Cuts & 
    \it Time & \multicolumn{1}{c}{\it Gap } &  \it Cuts  \\
    \hline
    6x12x1   & 0.60 & 0.00 & 1772 & 10 & 0.46 & 0.00 & 1619 & 10 & 0.05 & 7.79 & 803 & 0.05 & 5.13 & 809 \\
    6x12x2   & 0.45 & 0.00 & 1964 & 10 & 0.32 & 0.00 & 1861 & 10 & 0.07 & 4.28 & 1082 & 0.07 & 0.44 & 954 \\
    6x12x3   & 7254.68 & 2.49 & 33356 & 4 & 4449.73 & 0.73 & 14565 & 7 & 0.23 & 92.76 & 2490 & 0.14 & 14.10 & 1302 \\
    12x24x1   & 3.68 & 0.00 & 10745 & 10 & 3.30 & 0.00 & 10317 & 10 & 0.25 & 8.96 & 3330 & 0.23 & 8.33 & 3359 \\
    12x24x2   & 1148.31 & 0.16 & 23917 & 9 & 196.46 & 0.00 & 15677 & 10 & 0.32 & 8.34 & 3993 & 0.31 & 3.48 & 3782 \\
    12x24x3   & 10000.08 & 20.66 & 128509 & 0 & 10000.03 & 12.52 & 45311 & 0 & 1.26 & 99.80 & 5978 & 1.09 & 18.30 & 5062 \\
    24x48x1   & 123.30 & 0.00 & 65840 & 10 & 98.89 & 0.00 & 67111 & 10 & 1.91 & 15.04 & 15377 & 1.78 & 14.81 & 15401 \\
    24x48x2   & 10000.04 & 4.71 & 104814 & 0 & 10000.04 & 3.03 & 83944 & 0 & 1.95 & 12.29 & 15345 & 1.84 & 6.85 & 14754 \\
    24x48x3   & 10000.13 & 28.94 & 183575 & 0 & 10000.06 & 14.10 & 99062 & 0 & 10.38 & 99.92 & 20899 & 17.17 & 18.23 & 22284 \\
    \hline
    \end{tabular*}
}
\end{table}

The results show that MCM is always better than IM, even in the case of two segments: this is because, as Figure \ref{fig:ufl_1} shows, the function here is convex-concave rather than concave-convex and therefore does not satisfy the hypotheses of Proposition \ref{prop:two-segments}. Yet, the most significant improvements in the root node gap (from almost 100\% to around 10\%) are obtained when the number of segments is larger. This does not always result in a decreased running time due to the difficulty of the instances as the size grows, but at least it results in significantly lower final gaps.

\section{Conclusions}

Motivated by the SC-MINLP technique we have studied the tightness of the continuous relaxation bounds of the three ``textbook'' formulations for piecewise-convex functions. Unlike in the linear case, where they are all equivalent, one of them is shown to be weaker than the other two (unless in the concave-convex case with just two intervals), which does have a significant impact on actual computational results. We therefore believe that the result is significant, in that, given the well-known equivalence for the linear case, there may be a tendency to assume that the choice among these formulations is irrelevant: in fact, we precisely fell in this trap ourselves. Instead, the choice of the formulation is (as it often happens) important for final performances, and this may impact the many applications, besides SC-MINLP, where piecewise-convex functions are used.

\bibliographystyle{plain} 
\bibliography{main}

\appendix
\section*{Appendix}
Here we present the proof of Proposition~\ref{prop:two-segments}. As in the pictures, we consider the solution of the continuous relaxation of the two models as a function of $x$. To avoid having variables with the same names, for IM we replace $z^s$ with $\gamma^s$, $y^s$ with $\psi^s$, and $x^s$ with $\phi^s$. We, thus, have
\begin{align}%
g_{IM}(x) = \min \; & g(l^1) +  \alpha^1 \phi^1  + \gamma^2\\
& \gamma^2 \geq [ g(l^2 + \phi^2 / \psi^2 ) - g(l^2) ] \psi^2 
  \label{eq:gamma2_definition}\\
  & x= l^1 + \phi^1 + \phi^2  \label{eq:x_equal_phi1_plus_phi2}\\
  & (l^{s+1} - l^s ) \psi^{s+1} \leq \phi^s \leq (l^{s+1} - l^{s}) \psi^s
  & s=1,2 \label{eq:bound_constraints_IM}\\
  &  \psi^1 = 1 \;,\;  \psi^3 = 0 \label{eq:incremental_limits} \\
g_{MCM}(x) = \min \; & g(0) y^2 + \alpha^1 x^1 + [ g(l^1) - \alpha^1 l^1 ] y^1  + z^2 \\
    & x = x^1 + x^2 \label{eq:x_equal_x1_plus_x2}\\
    &  z^2 \geq [ g(x^2/ y^2 ) - g(0)] y^2 \\
    &  l^s y^s \leq x^s \leq l^{s+1} y^s
    &  s=1,2 \label{eq:bound_constraints_MCM}\\
    &  y^1 + y^2 = 1 \label{eq:multiple_choice_constraint}\\
    &  y^1, y^2 \geq 0 \label{eq:MCM_sign}
\end{align}
Note that in both models
\begin{equation}
  \alpha^1 = \frac{g(l^2) - g(l^1)}{l^2 - l^1}.
  \label{eq:alpha1_definition}
\end{equation}
Given an optimal solution $(\phi,\psi,\gamma)$ for $g_{IM}(x)$, $(x,y,z)$ defined by
\begin{align}
 & y^s = \psi^s - \psi^{s+1} \;\;,\;\;
   x^s = \phi^s + l^s \psi^s - l^{s+1} \psi^{s+1}
 & s=1,2 \label{eq:solxy_from_IM_to_MCM} \\
 & z^2 = [ g(x^2/y^2) - g(0) ] y^2 & \label{eq:solz_from_IM_to_MCM}
\end{align}
is feasible for the MCM. Indeed, \eqref{eq:x_equal_x1_plus_x2} is satisfied by \eqref{eq:solxy_from_IM_to_MCM} and \eqref{eq:x_equal_phi1_plus_phi2}:
\[ 
 x^1 + x^2 = \phi^1 + l^1 - l^2 \psi^2 + \phi^2 + l^2 \psi^2 = l^1 + \phi^1 + \phi^2 = x.
\]
For $s=1,2$, constraints \eqref{eq:bound_constraints_MCM} are satisfied because of \eqref{eq:solxy_from_IM_to_MCM} and \eqref{eq:bound_constraints_IM}:
\[
  \begin{array}{l}
    l^s (\psi^s - \psi^{s+1}) \leq \phi^s + l^s \psi^s - l^{s+1} \psi^{s+1}\leq l^{s+1}\psi^s - l^{s+1}\phi^{s+1} \equiv \\
    \equiv  
    (l^{s+1} - l^s) \phi^{s+1} \leq \phi^s \leq  (l^{s+1} - l^2) \phi^s
    \equiv l^s y^s \leq x^s \leq l^{s+1} y^s
  \end{array}
\] 
The multiple choice constraint~\eqref{eq:multiple_choice_constraint} is satisfied because of~\eqref{eq:solxy_from_IM_to_MCM} and \eqref{eq:incremental_limits}: 
\[
 y^1 + y^2 = 1 -  \psi^2 + \psi^2 = 1 
\] 
Then, the variable $z^2$ is computed through~\eqref{eq:solz_from_IM_to_MCM}, \eqref{eq:solxy_from_IM_to_MCM}, and~\eqref{eq:incremental_limits}:
\begin{equation}
 z^2 = [ g(x^2/y^2) - g(0)] y^2 =  [g(l^2+ \phi^2 / \psi^2 ) - g(0) ] \psi^2.
 \label{eq:solz_recomputed}
\end{equation}
Finally, the objective value of $(x,y,z)$ in the MC model is
\[
    \begin{array}{l}
    g(0) y^2 + \alpha^1 x^1 + [ g(l^1) - \alpha^1 l^1 ] y^1 + z^2 = \\
   g(0) \psi^2 + \alpha^1 ( \phi^1 + l^1 - l^2 \psi^2) + [ g(l^1) - \alpha^1 l^1](1-\psi^2) + z^2 = \\
    g(0) \psi^2 + \alpha^1 \phi^1 + g(l^1)  - g(l^1)\psi^2 - \alpha^1 (l^2 - l^1) \psi^2 + g(l^2 + \phi^2/\psi^2 ) \psi^2 - g(0) \psi^2 = \\
    g(l^1) + \alpha^1 \phi^1 - [g(l^1) + \frac{g(l^2) - g(l^1)}{l^2-l^1}(l^2-l^1)] \psi^2 + g(l^2 + \phi^2 / \psi^2 ) \psi^2 =\\
    g(l^1) + \alpha^1 \phi^1 + [ g(l^2 + \phi^2/\psi^2) - g(l^2)] \psi^2 \leq g(l^1) + \alpha^1 \phi^1 + \gamma^2 = g_{IM}(x),
    \end{array}
\]
where the first equality comes from~\eqref{eq:solxy_from_IM_to_MCM} and~\eqref{eq:incremental_limits}, the second equality comes from~\eqref{eq:solz_recomputed}, the third equality comes from~\eqref{eq:alpha1_definition}, the fourth equality comes from simple algebra, and the final inequality comes from \eqref{eq:gamma2_definition}. All this proves that $g_{MCM}(x) \leq g_{IM}(x)$: since  Corollary~\ref{cor:IMvsMCM} implies $g_{IM}(x) \leq g_{MCM}(x)$, necessarily $g_{IM}(x) = g_{MCM}(x)$.

Example~\ref{ex:convex-concave} shows that Proposition~\ref{prop:two-segments} does not hold if we switch the segments, i.e., $g$ is convex in $[l^1,l^2]$ and concave in $[l^2,l^3]$. One may still wonder why the above argument cannot be repeated in this case, where we rather have
\begin{align}
 g_{IM}(x) = \min \; & g(l^1) + \gamma^1 + \alpha^2 \phi^2 \nonumber \\
 & \gamma^1 \geq [g(l^1 + \phi^1/\psi^1) - g(l^1)] \psi^1 = g(l^1 + \phi^1) - g(l^1)
   \label{eq:IM_reversed} \\
 & \eqref{eq:x_equal_phi1_plus_phi2} \;\;,\;\; \eqref{eq:bound_constraints_IM} \;\;,\;\;
   \eqref{eq:incremental_limits} \nonumber \\
 g_{MCM}(x) = \min\; & g(0) y^1 + z^1 + \alpha^2 x^2 + [ g(l^2) - \alpha^2 l^2] y^2 \nonumber \\
 & z^1 \geq [g(x^1/y^1) - g(0)] y^1 \label{eq:MCM_reversed} \\
 & \eqref{eq:x_equal_x1_plus_x2} \;\;,\;\; \eqref{eq:bound_constraints_MCM} \;\;,\;\;
   \eqref{eq:multiple_choice_constraint} \;\;,\;\; \eqref{eq:MCM_sign} \nonumber
\end{align}
Comparing \eqref{eq:IM_reversed} with \eqref{eq:MCM_reversed} we see that the Perspective Reformulation has no effect on IM, while it strengthens MCM; the crucial equality in \eqref{eq:IM_reversed} comes from $\psi^1 = 1$ in \eqref{eq:incremental_limits}, i.e., \eqref{eq:Inc-firstlast}. Not fixing the variable would make the IM relaxation even weaker, and thus the result would then a fortiori hold.

\end{document}